\documentclass[11pt]{amsart}
\newtheorem{theorem}{Theorem}[section]
\newtheorem{corollary}[theorem]{Corollary}
\newtheorem{lemma}[theorem]{Lemma}
\newtheorem{proposition}[theorem]{Proposition}
\newtheorem{example}[theorem]{Example}
\theoremstyle{definition}

\theoremstyle{remark}

\numberwithin{equation}{section}

\begin{document}

\title[Dichotomy of stable radial solutions outside a ball]
{Dichotomy of stable radial solutions of $-\Delta u=f(u)$ outside a ball}
\author{Salvador Villegas}
\thanks{The author has been supported by the MEC Spanish grant MTM2012-37960.}
\address{Departamento de An\'{a}lisis
Matem\'{a}tico, Universidad de Granada, 18071 Granada, Spain.}
\email{svillega@ugr.es}

\begin{abstract}
This paper is devoted to the study of stable radial solutions of $-\Delta u=f(u) \mbox{ in } \mathbb{R}^N\setminus B_1=\{ x\in \mathbb{R}^N : \vert x\vert\geq 1\}$, where $f\in C^1(\mathbb{R})$ and $N\geq 2$. We prove that such solutions are either large [in the sense that $\vert u(r)\vert \geq M r^{-N/2+\sqrt{N-1}+2}\ $, if $2\leq N\leq 9$; $\vert u(r)\vert \geq M \log (r)\ $, if $N=10$; $\vert u(r)-u_\infty \vert \geq M r^{\-N/2+\sqrt{N-1}+2}\ $, if $N\geq 11$; $\forall r\geq r_0$, for some $M>0$, $r_0\geq 1$] or small [in the sense that $\vert u(r)\vert \leq M\log (r)\ $, if $N=2$; $\vert u(r)-u_\infty \vert \leq M r^{\-N/2-\sqrt{N-1}+2}$;\, if $N\geq 3$; $\forall r\geq 2$, for some $M>0$], where $u_\infty=\lim_{r\rightarrow \infty}u(r)\in [-\infty,+\infty]$. These results can be applied to stable outside a compact set radial solutions of equations of the type $-\Delta u=g(u) \mbox{ in } \mathbb{R}^N$. We prove also the optimality of these results, by considering solutions of the form $u(r)=r^\alpha$ or $u(r)=\log (r)$, $\forall r\geq 1$, where $\alpha \in \mathbb{R} \setminus \{ 0\}$.
\end{abstract}

\maketitle
\section{Introduction and main results}

This paper deals with the stability of radial solutions of

\begin{equation}\label{mainequation}
-\Delta u=f(u) \ \ \mbox{ in } \mathbb{R}^N\setminus B_1,
\end{equation}

\noindent where $B_1$ is the open unit ball of $\mathbb{R}^N$, $N\geq 2$ and $f\in C^1(\mathbb{R})$. We consider classical solutions $u\in C^2(\mathbb{R}^N\setminus B_1)$.

A solution $u$ of (\ref{mainequation}) is called stable if

$$Q_u(v):=\int_{\mathbb{R}^N\setminus \overline{B_1}} \left( \vert \nabla v\vert^2-f'(u)v^2\right) \, dx\geq 0$$

\noindent for every $v\in C^1 (\mathbb{R}^N\setminus \overline{B_1})$ with compact support in $\mathbb{R}^N\setminus \overline{B_1}$. Note that the above expression is nothing but the second variation of the energy functional associated to (\ref{mainequation}) in a bounded domain $\Omega\subset\mathbb{R}^N\setminus \overline{B_1}$: $E_\Omega (u)=\int_\Omega \left( \vert \nabla u\vert^2 /2-F(u)\right) \, dx$, where $F'=f$. Thus, if $u\in C^1(\mathbb{R}^N\setminus B_1)$ is a local minimizer of $E_\Omega$ for every bounded smooth domain $\Omega\subset\mathbb{R}^N\setminus\overline{B_1}$ (i.e., a minimizer under every small enough $C^1(\overline{\Omega})$ perturbation vanishing on $\partial \Omega$), then $u$ is a stable solution of (\ref{mainequation}).

We will be also interested in stable outside a compact set radial solutions of
\begin{equation}\label{space}
-\Delta u=g(u) \ \ \mbox{ in } \mathbb{R}^N,
\end{equation}

\noindent where  $N\geq 2$ and $g\in C^1(\mathbb{R})$.

We say that a classical solution $u\in C^2(\mathbb{R}^N)$ of (\ref{space}) is stable outside a compact set if there exists a compact set $K\subset\mathbb{R}^N$ such that $Q_u(v)=\int_{\mathbb{R}^N}\left( \vert \nabla v\vert^2-g'(u)v^2\right) \, dx\geq 0$ for every $v\in C^1 (\mathbb{R}^N)$ with compact support in $\mathbb{R}^N\setminus K$.

Clearly the stability outside a compact set of a solution of (\ref{space}) is equivalent to the existence of $R_0>0$ such that $u$ is stable in $\mathbb{R}^N\setminus B_{R_0}$. It follows easily that the function $w(x):= u(R_0 x)$ is an stable solution of $-\Delta w=R_0^2 g(w)$ in $\mathbb{R}^N\setminus B_1$ and we can apply the results obtained for such solutions.

On the other hand we say that a classical solution $u\in C^2(\mathbb{R}^N)$ of (\ref{space}) has finite Morse index equal to an integer $k\geq 0$ if $k$ is the maximal dimension of a subspace $X_k\subset C_c^1(\mathbb{R}^N)$ (the space of  $C^1(\mathbb{R}^N)$ functions with compact support) such that

$$Q_u(\varphi)=\int_{\mathbb{R}^N} \left( \vert \nabla \varphi\vert^2-g'(u)\varphi^2\right) \, dx<0 \ \ \ \mbox{ for all }\varphi\in X_k\setminus\{ 0\}.$$

If there is no such finite integer $k$, we say that $u$ has infinite Morse index.

Clearly, every stable solution has finite Morse index equal to $0$. It is also easily seen that every solution with finite Morse index is stable outside a compact set. Indeed, if $X_k=\mbox{Span}\left\{ \varphi_1,\ldots, \varphi_k\right\}$ is a subspace of dimension $k$ of $C_c^1(\mathbb{R}^N)$ such that $Q_u(\varphi)<0$ for any $\varphi\in X_k\setminus\{ 0\}$ and $K:=\bigcup_{j=1}^k \mbox{supp}(\varphi_j)$, then $Q_u(v)\geq 0$ for every $v\in C_c^1(\mathbb{R}^N\setminus K)$, and the claim is proved. Hence, we can apply to finite Morse index solutions the result obtained for solutions which are stable outside a compact set.

Farina \cite{farina1, farina2} studied the stability and stability outside a compact set of nontrivial solutions of the Lane-Emden equation $-\Delta u=\vert u\vert^{p-1}u$ in $\mathbb{R}^N$ ($p>1$). It is proved that the existence of such solutions depend on $N$ (dimensions $N=2$ and $N=10$ are critical in some sense, as in the main results of this paper) and $p$ (there are two critical values: $p=(N+2)/(N-2)$, the usual critical exponent in Sobolev imbedding theorems, defined for $N\geq 3$; and $p=p_c:=\left((N-2)^2-4N+8\sqrt{N-1}\right)/((N-2)(N-10))$, defined for $N\geq 11$). A complete classification of radial solutions of this equation which are stable outside a compact set is also given (see Remark 4 below).

Farina \cite{farina3} considered the equation $-\Delta u=e^u$ in $\mathbb{R}^N$, obtaining that there are no stable solutions if $N\leq 9$. A complete classification of solutions of this equation which are stable outside a compact set is obtained for $N=2$ (these solutions are radially symmetric, up to a translation). In a later paper, Dancer and Farina \cite{dancerfarina} studied also this equation and proved that there are no solutions which are stable outside a compact set if $3\leq N\leq 9$ (again dimensions $N=2$ and $N=10$ are critical in some sense).

Dupaigne and Farina \cite{dupaignefarina1, dupaignefarina2} have also studied the stability and stability outside a compact set of equation (\ref{space}), for a large class of functions $g\in C^1(\mathbb{R})$. Among other things, they proved that if $g\geq 0$ and $1\leq N\leq 4$, then there are no nonconstant bounded stable solutions.

We are interested in radial solutions of (\ref{mainequation}) and (\ref{space}). By abuse of notation, we write $u(r)$ instead of $u(x)$, where $r=\vert x\vert$ and $x\in \mathbb{R}^N$. We will denote by $u_r$ the radial derivative of a radial function u.

Concerning with nonconstant bounded stable radial solutions $u\in C^2(\mathbb{R}^N)$ of (\ref{space}), Cabr\'e and Capella \cite{cabrecapella1} proved that there are no such solutions if $N\leq 10$ and $g$ satisfies a nondegeneracy condition. The author \cite{yo} refined this result, proving that this nondegeneracy condition is not necessary and giving sharp pointwise estimates related to the asymptotic behavior of such solutions (not necessarily bounded). Specifically, in \cite{yo} it is proved that every nonconstant radial stable solution of (\ref{space}) satisfies $\vert u(r)\vert \geq M r^{-N/2+\sqrt{N-1}+2}\ $, if $N\neq 10$, and $\vert u(r)\vert \geq M \log (r)\ $, if $N=10$; $\forall r\geq r_0$, for some $M,r_0>0$.

In this paper we establish that there is a dichotomy of radial stable solutions of (\ref{mainequation}): such solutions are either {\it large} (i.e. roughly speaking, they grow at least like a power $r^{-N/2+\sqrt{N-1}+2}$, like the mentioned stable solutions in all of $\mathbb{R}^N$) or {\it small} (i.e. roughly speaking, they grow at most like a power $r^{-N/2-\sqrt{N-1}+2}$). Note that the exponents $-N/2\pm\sqrt{N-1}+2$ vanish at $N=2$ and at $N=10$, respectively. Hence, we can consider that these dimensions are critical in the study of these problems.

\begin{theorem}\label{dichotomy}

Let $N\geq 2$, $f\in C^1(\mathbb{R})$ and $u$ be a radial stable solution of (\ref{mainequation}) (not necessarily bounded). Then there exists $u_\infty:=\lim_{r\to \infty}u(r)\in [-\infty,+\infty]$ and $u$ satisfies either $(L)$ or $(S)$. Here

\

$(L)$ There exist $M>0$, $r_0\geq 1$ such that
$$
\begin{array}{rll}
\vert u(r)\vert & \geq M r^{-N/2+\sqrt{N-1}+2}& \forall r\geq r_0 \mbox{ if } 2\leq N\leq 9,\\
\vert u(r)\vert & \geq M \log (r) & \forall r\geq r_0 \mbox{ if } N=10,\\
\vert u(r)-u_\infty\vert & \geq M r^{-N/2+\sqrt{N-1}+2}& \forall r\geq r_0 \mbox{ if } N\geq 11.\\
\end{array}
$$

\

$(S)$ There exists $M>0$ such that
$$
\begin{array}{rll}
\vert u(r)\vert & \leq M \log (r) & \forall r\geq 2 \mbox{ if } N=2, \\
\vert u(r)-u_\infty\vert & \leq M r^{-N/2-\sqrt{N-1}+2}& \forall r\geq 1\mbox{ if } N\geq 3.\\
\end{array}
$$

\end{theorem}

\

Note that in the case $N\geq 11$, condition $(L)$ is relevant if $u$ is bounded. Otherwise $u_\infty=\pm \infty$ and the inequality in $(L)$ is vacuous. On the other hand, if $N\geq 3$, condition $(S)$ says implicitly that $u_\infty\in\mathbb{R}$ and hence $u$ is bounded.

Example \ref{optimalexample} below shows that the exponents $-N/2\pm \sqrt{N-1}+2$ which appear in Theorem \ref{dichotomy} are optimal. In fact any pure power in the set $(-\infty, -N/2-\sqrt{N-1}+2] \cup [-N/2+\sqrt{N-1}+2,+\infty)$ is allowed for stable solutions (considering logarithm as a $0$-power).

If $u\in C^2(\mathbb{R}^N$) is a nonconstant radial stable solution of (\ref{space}) then, by a result of the author (\cite[Th. 2.1 and 2.2]{yo}), $u$ satisfies $(L)$. This is consistent with Theorem \ref{dichotomy}.

\begin{theorem}\label{dichotomybounded}

Let $N\geq 3$, $f\in C^1(\mathbb{R})$ and $u$ be a bounded radial stable solution of (\ref{mainequation}). Then there exists $u_\infty:=\lim_{r\to \infty}u(r)\in\mathbb{R}$ and $u$ satisfies either $(L')$ or $(S')$. Here

\

$(L')$ $N\geq 11$ and there exist $M>0$, $r_0\geq 1$ such that
$$
\vert u(r)-u_\infty\vert\geq M r^{-N/2+\sqrt{N-1}+2} \ \ \forall r\geq r_0.
$$

\ \ \ \ \ \ Moreover, $\vert\nabla u\vert\notin L^p(\mathbb{R}^N\setminus\overline{B_1})$ for every $1\leq p\leq \frac{N}{N/2-\sqrt{N-1}-1}$.

\ \ \ \ \ \ In particular $\vert\nabla u\vert\notin L^2(\mathbb{R}^N\setminus\overline{B_1})$ .

\

$(S')$ There exists $M>0$ such that
$$
\begin{array}{ll}
\vert u(r)-u_\infty\vert\leq M r^{-N/2-\sqrt{N-1}+2} & \forall r\geq 1.\\
\end{array}
$$

\ \ \ \ \ \ Moreover, $\vert\nabla u\vert\in L^p(\mathbb{R}^N\setminus \overline{B_1})$ for every \ $\frac{N}{N/2+\sqrt{N-1}-1}<p\leq \infty$.

\ \ \ \ \ \ In particular $\vert\nabla u\vert\in L^2(\mathbb{R}^N\setminus\overline{B_1})$.

\end{theorem}

\

\begin{corollary}\label{morse}

Let $N\geq 2$, $g\in C^1(\mathbb{R})$ and $u$ be a radial solution of (\ref{space}) which is stable outside a compact set (not necessarily bounded). Then $u$ satisfies either $(L)$ or $(S)$.

\end{corollary}

\begin{corollary}\label{morsebounded}

Let $N\geq 3$, $g\in C^1(\mathbb{R})$ and $u$ be a bounded radial solution of (\ref{space}) which is stable outside a compact set. Then $u$ satisfies either $(L')$ or $(S')$.

\end{corollary}

To distinguish if a nonconstant radial stable solution of (\ref{mainequation}) is {\it large} or {\it small} we will consider the following properties:

\

\noindent $(H_L)$ \ \ There exists $R_1\geq 1$ such that $$\int_{R_1}^{R_2} r^{N-1}u_r^2 \left( \eta'^2-\frac{N-1}{r^2}\eta^2\right) dr \geq 0,$$

\ \ \ \ \ \ \ for every  $R_2>R_1$ and $\eta\in C^{0,1}([R_1,R_2])$ such that $\eta(R_2)=0$.

\

\noindent $(H_S)$ \ \ For every $R_1\geq 1$ there exist $R_2>R_1$ and $\eta_0\in C^{0,1}([R_1,R_2])$

\ \ \ \ \ \ such that $\eta_0(R_2)=0$ and $$\int_{R_1}^{R_2} r^{N-1}u_r^2 \left( \eta_0'^2-\frac{N-1}{r^2}\eta_0^2\right) dr<0.$$

\

Note that $(H_L)$ and $(H_S)$ are complementary properties. We will show that a nonconstant radial stable solution of (\ref{mainequation}) satisfying $(H_L)$ is {\it large} (i.e. satisfies $(L)$) while a nonconstant radial stable solution of (\ref{mainequation}) satisfying $(H_S)$ is {\it small} (i.e. satisfies $(S)$).

\section{Proof of the main results}

The following lemma follows easily from the ideas of the proof of \cite[Lem. 2.2]{cabrecapella1}, which was inspired by the proof of Simons theorem on the nonexistence of singular minimal cones in $\mathbb{R}^N$ for $N\leq 7$ (see \cite[Th. 10.10]{simons} and \cite[Rem. 2.2]{cabrecapella2} for more details).

\begin{lemma}\label{maintool}

Let $N\geq 2$, $f\in C^1(\mathbb{R})$ and $u$ be a radial stable solution of (\ref{mainequation}).  Let $1\leq r_1<r_2<\infty$ and $\eta\in C^{0,1}([r_1,r_2])$ such that $\eta u_r$ vanishes at $r=r_1$ and $r=r_2$. Then

$$\int_{r_1}^{r_2} r^{N-1}u_r^2 \left( \eta'^2-\frac{N-1}{r^2}\eta^2\right) dr\geq 0.$$

\end{lemma}

\noindent\textbf{Proof.}
First of all, note that we can extend the second variation of energy $Q_u$ to the set of functions $v\in C^{0,1}(\mathbb{R}^N\setminus\overline{B_1})$ with compact support in $\mathbb{R}^N\setminus\overline{B_1}$, obtaining $Q_u(v)\geq 0$ for such functions $v$. Hence, if $r_1>1$, we can take the radial function $v=\eta u_r \chi_{B_{r_2}\setminus\overline{B_{r_1}}}$. In fact, by an approximative method, we can also take this function $v$ in the case $r_1=1$.

On the other hand, differentiating (\ref{mainequation}) with respect to r, we have

$$-\Delta u_r+\frac{N-1}{r^2}u_r=f'(u)u_r, \ \ \mbox{ for all }r\geq 1.$$

Following the ideas of the proof of \cite[Lem. 2.2]{cabrecapella1}, we can multiply this equality by $\eta^2 u_r$ and integrate by parts in the annulus of radii $r_1$ and $r_2$ to obtain

$$0=\int_{B_{r_2}\setminus\overline{B_{r_1}}}\left(\nabla u_r\nabla\left(\eta^2 u_r\right)+\frac{N-1}{r^2}u_r \eta^2 u_r-f'(u)u_r\eta^2 u_r\right)dx$$
$$=\int_{B_{r_2}\setminus\overline{B_{r_1}}}\left(\vert \nabla\left(\eta u_r\right)\vert^2-f'(u)\left( \eta u_r\right)^2\right)dx-\int_{B_{r_2}\setminus\overline{B_{r_1}}} u_r^2\left(\vert\nabla\eta\vert^2-\frac{N-1}{r^2}\eta^2\right)dx$$
$$=Q(\eta u_r \chi_{B_{r_2}\setminus\overline{B_{r_1}}})-\omega_N\int_{r_1}^{r_2} r^{N-1}u_r^2 \left( \eta'^2-\frac{N-1}{r^2}\eta^2\right) dr.$$

Using the stability of $u$ the lemma follows. \qed

\begin{lemma}\label{monotonicity}

Let $N\geq 2$, $f\in C^1(\mathbb{R})$ and $u$ be a nonconstant radial stable solution of (\ref{mainequation}). Then $u_r$ vanishes at most in one value in $[1,+\infty)$.

\end{lemma}

\noindent\textbf{Proof.}
Suppose by contradiction that there exist $1\leq r_1<r_2<\infty$ such that $u_r(r_1)=u_r(r_2)=0$. Taking $\eta\equiv 1$ in the previous lemma, we obtain

$$\int_{r_1}^{r_2} r^{N-1}u_r^2 \left( -\, \frac{N-1}{r^2}\right) dr\geq 0.$$

Hence we conclude that $u_r\equiv 0$ in $[r_1,r_2]$, which clearly forces $u$ is constant in $\mathbb{R}^N\setminus B_1$, a contradiction. \qed

\subsection{Large solutions}

In this subsection we will prove that a nonconstant radial stable solution of (\ref{mainequation}) satisfying $(H_L)$ is {\it large} (i.e. satisfies $(L)$).

\begin{lemma}\label{largelemma1}
Let $N\geq 2$ and $u$ be a nonconstant radial stable solution of (\ref{mainequation}) satisfying $(H_L)$. Then there exist $a\geq 1$ and $K>0$ such that

$$\int_r^{2r}\frac{ds}{u_r(s)^2}\leq K r^{N-2\sqrt{N-1}-1}\, \ \ \ \forall r\geq a.$$

\end{lemma}

\

\noindent\textbf{Proof.}
Consider $R_1\geq 1$ of $(H_L)$. From Lemma \ref{monotonicity} we can choose $a>R_1$ such that $u_r$ does not vanish in $[a,\infty)$. We now fix $r\geq a$ and consider the function

$$\eta (t)=\left\{
\begin{array}{ll}
a^{-\sqrt{N-1}} & \mbox{ if } R_1\leq t <a, \\ \\
t^{-\sqrt{N-1}} & \mbox{ if } a\leq t<r, \\ \\
\displaystyle{\frac{r^{-\sqrt{N-1}}}{\int_r^{2r}
\frac{ds}{u_r(s)^2}}} \int_t^{2r} \frac{ds}{u_r(s)^2} &
\mbox{ if } r\leq t\leq 2r. \\
\end{array}
\right.
$$

Applying $(H_L)$ (with $R_2=2r$) we have

$$0\leq \int_{R_1}^{R_2} t^{N-1}u_r(t)^2 \left( \eta'(t)^2-\frac{N-1}{t^2}\eta(t)^2\right) dt$$
$$=-(N-1)a^{-2\sqrt{N-1}}\int_{R_1}^{a}t^{N-3}u_r(t)^2 dt+\int_r^{2r} t^{N-1}u_r(t)^2 \left( \eta'(t)^2-\frac{N-1}{t^2}\eta(t)^2\right) dt$$
$$\leq-(N-1)a^{-2\sqrt{N-1}}\int_{R_1}^{a}t^{N-3}u_r(t)^2 dt+(2r)^{N-1}\int_r^{2r}u_r(t)^2 \eta'(t)^2 dt$$
$$=-(N-1)a^{-2\sqrt{N-1}}\int_{R_1}^{a}t^{N-3}u_r(t)^2 dt+(2r)^{N-1}\displaystyle{\frac{r^{-2\sqrt{N-1}}}{\int_r^{2r}\frac{ds}{u_r(s)^2}}}.$$

This gives

$$(N-1)a^{-2\sqrt{N-1}}\int_{R_1}^{a}t^{N-3}u_r(t)^2 dt\leq 2^{N-1}\displaystyle{\frac{r^{N-2\sqrt{N-1}-1}}{\int_r^{2r}\frac{ds}{u_r(s)^2}}},$$

\noindent which is the desired conclusion for

$$K=2^{N-1}/\left((N-1)a^{-2\sqrt{N-1}}\int_{R_1}^{a}t^{N-3}u_r(t)^2 dt\right).$$ \qed

\begin{lemma}\label{largelemma2}

Let $N\geq 2$ and $u$ be a nonconstant radial stable solution of (\ref{mainequation}) satisfying $(H_L)$. Then there exist $a\geq 1$ and $M'>0$ such that

$$\vert u(2r)-u(r)\vert \geq M' r^{-N/2+\sqrt{N-1}+2}\,  \ \ \
\forall r\geq a.$$
\end{lemma}

\

\noindent\textbf{Proof.}
Take the same constant $a\geq 1$ of Lemma \ref{largelemma1}. Fix $r\geq a$ and consider the functions:

$$\begin{array}{ll}
 \alpha (s)=\vert u_r(s)\vert^{-\frac{2}{3}}, & s\in (r,2r). \\
 \\
  \beta (s)=\vert u_r(s)\vert^{\frac{2}{3}}, & s\in (r,2r). \\

\end{array}
$$

By Lemma \ref{largelemma1} we have

$$\Vert \alpha \Vert _{L^3(r,2r)}\leq K^\frac{1}{3} \,
r^\frac{N-2\sqrt{N-1}-1}{3}$$

\noindent for a constant $K>0$ not depending on $r\geq a$. On the other hand, since $u_r$ does not vanish in $[a, \infty )$, it follows

$$\Vert \beta \Vert _{L^{3/2}(r,2r)}=\vert u(2r)-u(r)
\vert^\frac{2}{3}.$$

Applying H{\"o}lder inequality to functions $\alpha$ and $\beta$ we deduce

$$r=\int_r^{2r} \alpha(s)\beta(s) ds\leq \Vert \alpha \Vert _{L^3(r,2r)} \Vert \beta \Vert _{L^{3/2}(r,2r)}\leq K^\frac{1}{3}r^\frac{N-2\sqrt{N-1}-1}{3}\, \vert u(2r)-u(r)
\vert^\frac{2}{3},$$

\noindent which is the desired conclusion for $M'=K^{-1/2}$. \qed

\begin{proposition}\label{largeproposition}

Let $N\geq 2$ and $u$ be a nonconstant radial stable solution of (\ref{mainequation}) satisfying $(H_L)$. Then $u$ satisfies $(L)$.

\end{proposition}

\noindent\textbf{Proof.}
Consider the numbers $a\geq 1$ and $M'>0$ of Lemma \ref{largelemma2}. The proof will be divided into three cases:

\begin{enumerate}

\item[$\bullet$] Case $2\leq N\leq 9$.

It is easily seen that for every $r\geq a$ there exist an integer $m\geq 0$ and $a\leq z<2a$ such that $r=2^mz$. Thus, from Lemma \ref{largelemma2} and the monotonicity of $u$ in $[a,\infty)$, it follows that

$$\vert u(r)\vert \geq \vert u(r)-u(z)\vert -\vert u(z)\vert
=\sum_{k=1}^m\vert u(2^kz)-u(2^{k-1}z)\vert \, -\vert
u(z)\vert $$

$$\geq M'\sum_{k=1}^m(2^{k-1}z)^{-N/2+\sqrt{N-1}+2}\, -\vert u(z)\vert$$

$$=M'\left(\frac{
r^{-N/2+\sqrt{N-1}+2}-z^{-N/2+\sqrt{N-1}+2}}{2^{-N/2+\sqrt{N-1}+2}-1}\right)
-\vert u(z)\vert \, ,$$

\noindent where $M'>0$ does not depend on $r\geq a$. Since $z\in [a,2a)$, $u$ is continuous and $-N/2+\sqrt{N-1}+2>0$, the above inequality is of the type

$$\vert u(r)\vert \geq M_1 r^{-N/2+\sqrt{N-1}+2}-M_2 \, \ \ \ \ \
\ \forall r\geq a,$$

\noindent for certain $M_1, M_2>0$. It follows easily $(L)$ in this case.

\item[$\bullet$] Case $N=10$.

In this case $-N/2+\sqrt{N-1}+2=0$. Following the same notation of the previous case, we can apply the some reasoning and conclude that

$$\vert u(r)\vert \geq M'm-\vert u(z)\vert =\frac{M'(\log r -\log z)}{\log 2}-\vert u(z)\vert ,$$

\noindent and $(L)$ follows immediately for this case.

\item[$\bullet$] Case $N\geq 11$.

From Lemma \ref{monotonicity} we deduce that there exists $u_\infty:=\lim_{r\to \infty}u(r)\in [-\infty,+\infty]$. If $u_\infty=\pm \infty$ then the inequality in $(L)$ is trivial. Then without loss of generality we can assume that $u_\infty\in\mathbb{R}$. Let $r\geq a$. From Lemma \ref{largelemma2} and the monotonicity of $u$ in $[a,\infty)$ we see that

$$\vert u_\infty -u(r)\vert =\sum_{k=1}^{\infty}\vert u(2^k
r)-u(2^{k-1}r)\vert\geq
M'\sum_{k=1}^{\infty}(2^{k-1}r)^{-N/2+\sqrt{N-1}+2}$$

$$=\left( M'\sum_{k=1}^{\infty}2^{(k-1)(-N/2+\sqrt{N-1}+2)}\right) r^{-N/2+\sqrt{N-1}+2}.$$

Finally, since $-N/2+\sqrt{N-1}+2<0$, the above series is convergent and $(L)$ is proved in this case with $r_0=a$. \qed

\end{enumerate}

\subsection{Small solutions}

In this subsection we will prove that a radial stable solution of (\ref{mainequation}) satisfying $(H_S)$ is {\it small} (i.e. satisfies $(S)$).

\begin{lemma}\label{smalllemma1}
Let $N\geq 2$ and $u$ be a radial stable solution of (\ref{mainequation}) satisfying $(H_S)$. Then there exists $K>0$ such that

$$\int_r^{2r}u_r(s)^2 ds \leq K r^{-N-2\sqrt{N-1}+3}\, \ \ \ \forall r\geq 1.$$
\end{lemma}

\

\noindent\textbf{Proof.}
Take an arbitrary $r\geq 2$ and consider $R_1=2r$ in $(H_S)$. Then there exist $R_2>2r$ and $\eta_0\in C^{0,1}([2r,R_2])$ such that $\eta_0(R_2)=0$ and

\begin{equation}\label{yyy}\int_{2r}^{R_2} t^{N-1}u_r^2(t) \left( \eta_0'(t)^2-\frac{N-1}{t^2}\eta_0(t)^2\right) dt<0.\end{equation}

Note that $\eta_0(2r)\neq 0$ (otherwise we would obtain a contradiction with Lemma \ref{maintool} for $\eta=\eta_0$ and $[r_1,r_2]=[2r,R_2]$). Thus, multiplying by a constant if necessary, there is no loss of generality in assuming $\eta_0(2r)=r^{\sqrt{N-1}}$. We now fix $r\geq 2$ and consider the function

$$\eta (t)=\left\{
\begin{array}{ll}
2^{\sqrt{N-1}}(t-1) & \mbox{ if } 1\leq t <2, \\ \\
t^{\sqrt{N-1}} & \mbox{ if } 2\leq t<r, \\ \\
r^{\sqrt{N-1}} & \mbox{ if } r\leq t<2r, \\ \\
\eta_0(t) & \mbox{ if } 2r\leq t\leq R_2.
\end{array}
\right.
$$

Applying (\ref{yyy}) and Lemma \ref{maintool} to this function $\eta\in C^{0,1}([1,R_2])$ we have

$$0\leq \int_1^{R_2} t^{N-1}u_r(t)^2 \left( \eta'(t)^2-\frac{N-1}{t^2}\eta(t)^2\right) dt$$
$$=4^{\sqrt{N-1}}\int_1^2 t^{N-1}u_r(t)^2\left(1-\frac{(N-1)(t-1)^2}{t^2}\right)dt-(N-1)r^{2\sqrt{N-1}}\int_r^{2r} t^{N-3}u_r(t)^2 dt$$ $$+\int_{2r}^{R_2} t^{N-1}u_r(t)^2 \left( \eta_0'(t)^2-\frac{N-1}{t^2}\eta_0(t)^2\right) dt$$ $$<4^{\sqrt{N-1}}\int_1^2 t^{N-1}u_r(t)^2\left(1-\frac{(N-1)(t-1)^2}{t^2}\right)dt-(N-1)r^{2\sqrt{N-1}}C_N r^{N-3}\int_r^{2r}u_r(t)^2 dt,$$

\noindent where $C_N=\min\{1, 2^{N-3}\}$. This gives

$$(N-1)r^{2\sqrt{N-1}}C_N r^{N-3}\int_r^{2r}u_r(t)^2 dt<4^{\sqrt{N-1}}\int_1^2 t^{N-1}u_r(t)^2\left(1-\frac{(N-1)(t-1)^2}{t^2}\right)dt,$$

\noindent which is our claim (if $r\geq 2$) for

$$K=4^{\sqrt{N-1}}\int_1^2 t^{N-1}u_r(t)^2\left(1-\frac{(N-1)(t-1)^2}{t^2}\right)dt/\left((N-1)C_N\right).$$

Finally, if $1\leq r<2$, since $u_r$ is bounded in the interval $[1,4]$, we also have the desired inequality and the lemma follows easily. \qed

\begin{lemma}\label{smalllemma2}

Let $N\geq 2$ and $u$ be a radial stable solution of (\ref{mainequation}) satisfying $(H_S)$. Then there exists $M'>0$ such that

$$\vert u(2r)-u(r)\vert \leq M' r^{-N/2-\sqrt{N-1}+2}\,  \ \ \
\forall r\geq 1.$$
\end{lemma}

\

\noindent\textbf{Proof.}
Fix $r\geq 1$. Applying Lemma \ref{smalllemma1} and Cauchy-Schwarz inequality we obtain

$$\vert u(2r)-u(r)\vert=\left\vert \int_r^{2r}u_r(s)ds\right\vert\leq\int_r^{2r}\vert u_r(s)\vert ds\leq \left( \int_r^{2r}u_r(s)^2ds\right)^{1/2}\left( \int_r^{2r}ds\right)^{1/2}$$

$$\leq K^{\frac{1}{2}} r^{\frac{-N-2\sqrt{N-1}+3}{2}}r^{\frac{1}{2}},$$

\noindent which is our assertion for $M'=K^{1/2}$. \qed

\begin{proposition}\label{smallproposition}

Let $N\geq 2$ and $u$ be a radial stable solution of (\ref{mainequation}) satisfying $(H_S)$. Then $u$ satisfies $(S)$.

\end{proposition}

\noindent\textbf{Proof.}
The proof will be divided into two cases:

\begin{enumerate}

\item[$\bullet$] Case $N=2$.

In this case $-N/2-\sqrt{N-1}+2=0$. Let $r\geq 1$. Then there exist an integer $m\geq 0$ and $1\leq z<2$ such that $r=2^mz$. Thus, from Lemma \ref{smalllemma2} it follows that

$$\vert u(r)\vert \leq \vert u(r)-u(z)\vert +\vert u(z)\vert \leq \sum_{k=1}^m\vert u(2^kz)-u(2^{k-1}z)\vert \, +\vert
u(z)\vert\leq M'm+\vert u(z)\vert$$ $$ =\frac{M'(\log r -\log z)}{\log 2}+\vert u(z)\vert ,$$

\noindent where $M'>0$ does not depends on $r\geq 1$. Since $z\in [1,2)$ and $u$ is continuous, then $(S)$ follows immediately for this case.

\item[$\bullet$] Case $N\geq 3$.

Let $r\geq 1$ and $j\in \mathbb{N}$ be arbitrary. Using Lemma \ref{smalllemma2} we obtain

$$\vert u(2^j r)-u(r)\vert\leq\sum_{k=1}^j\vert u(2^k r)-u(2^{k-1}r)\vert\leq M'\sum_{k=1}^j (2^{k-1}r)^{-N/2-\sqrt{N-1}+2}$$ $$=\left(M'\sum_{k=1}^j 2^{(k-1)(-N/2-\sqrt{N-1}+2)}\right) r^{-N/2-\sqrt{N-1}+2}.$$

Since $-N/2-\sqrt{N-1}+2<0$ the above series is convergent and we can let $j\to \infty$ (remember that the existence of $u_\infty:=\lim_{r\to \infty}u(r)\in [-\infty,+\infty]$ is guaranteed by Lemma \ref{monotonicity}) to obtain

$$\vert u_\infty-u(r)\vert\leq \frac{M'}{1-2^{-N/2-\sqrt{N-1}+2}}r^{-N/2-\sqrt{N-1}+2},$$

\noindent and the proof is complete in this case. \qed

\end{enumerate}

\subsection{Proofs of Theorems \ref{dichotomy} and \ref{dichotomybounded} and Corollaries \ref{morse} and \ref{morsebounded}}
In this subsection we use the obtained results for {\it large} and {\it small} solutions to prove the main results of this paper.

\

\noindent\textbf{Proof of Theorem \ref{dichotomy}.}
Let $N\geq 2$, $f\in C^1(\mathbb{R})$ and $u$ be a radial stable solution of (\ref{mainequation}) (not necessarily bounded). By Lemma \ref{monotonicity} it is immediate that there exists $u_\infty:=\lim_{r\to \infty}u(r)\in [-\infty,+\infty]$.

If $u$ is constant, then clearly $u$ satisfies $(S)$ and not $(L)$. Hence, from the rest of the proof we will suppose that $u$ is not constant. Obviously $(H_L)$ and $(H_S)$ are complementary properties, i.e. $u$ satisfies either $(H_L)$ or $(H_S)$. By Proposition \ref{largeproposition} if $u$ satisfies $(H_L)$ then $u$ satisfies $(L)$ and by Proposition \ref{smallproposition} if $u$ satisfies $(H_S)$ then $u$ satisfies $(S)$.

To finish the proof, let us observe that conditions $(L)$ and $(S)$ are clearly incompatible. \qed

\

\noindent\textbf{Proof of Theorem \ref{dichotomybounded}.}
Let $N\geq 3$, $f\in C^1(\mathbb{R})$ and $u$ be a bounded radial stable solution of (\ref{mainequation}). By Theorem \ref{dichotomy} we have that there exists $u_\infty:=\lim_{r\to \infty}u(r)\in \mathbb{R}$.  If $u$ is constant, then clearly $u$ satisfies $(S')$ and not $(L')$. Hence, from the rest of the proof we will suppose that $u$ is not constant.

Suppose that $u$ satisfies $(H_L)$. We will show that $u$ satisfies $(L')$. Since $u$ satisfies $(L)$ (again from Proposition \ref{largeproposition}) and $-N/2+\sqrt{N-1}+2>0$ for $3\leq N\leq 9$, we deduce $N\geq 11$. What is left to show in this case is that $\vert\nabla u\vert\notin L^p(\mathbb{R}^N\setminus\overline{B_1})$ for every $1\leq p\leq N/\left(N/2-\sqrt{N-1}-1\right)$. To this end, take the constant $a\geq 1$ of Lemma \ref{largelemma1}, fix $r\geq a$ and consider the functions

$$\begin{array}{ll}
 \alpha (s)=\vert u_r(s)\vert^{-\frac{2p}{p+2}}, & s\in (r,2r). \\
 \\
  \beta (s)=\vert u_r(s)\vert^{\frac{2p}{p+2}}, & s\in (r,2r). \\
\end{array}$$

Applying Lemma \ref{largelemma1} and H{\"o}lder inequality to functions $\alpha$ and $\beta$ we deduce

$$r=\int_r^{2r} \alpha(s)\beta(s) ds\leq \Vert \alpha \Vert_{L^{\frac{p+2}{p}}(r,2r)} \Vert \beta \Vert_{L^\frac{p+2}{2}(r,2r)}$$
$$\leq K^\frac{p}{p+2} \, r^{\frac{p}{p+2}\left(N-2\sqrt{N-1}-1\right)}\left(\int_r^{2r} \vert u_r(s)\vert^p ds\right)^{\frac{2}{p+2}},$$

\noindent for a constant $K>0$ not depending on $r\geq a$.

This gives

$$\int_r^{2r} \vert u_r(s)\vert^p ds\geq K^{\frac{-p}{2}}r^{\, p\left(-N/2+\sqrt{N-1}+1\right)+1},$$

\noindent for $r\geq a$. From this, we obtain

$$\int_r^{2r} s^{N-1} \vert u_r(s)\vert^p ds\geq r^{N-1} \int_r^{2r} \vert u_r(s)\vert^p ds\geq K^{\frac{-p}{2}}r^{\, p\left(-N/2+\sqrt{N-1}+1\right)+N},$$

\noindent for $r\geq a$. Finally, since $N\geq 11$, $1\leq p\leq N/\left(N/2-\sqrt{N-1}-1\right)$ and $r\geq a\geq 1$ we deduce that $r^{\, p\left(-N/2+\sqrt{N-1}+1\right)+N}\geq 1$, which implies
$\int_r^{2r} s^{N-1} \vert u_r(s)\vert^p ds\geq K^{-p/2}$ for a constant $K>0$ not depending on $r\geq a$. We conclude that $s^{N-1}\vert u_r(s)\vert^p\notin L^1(a,\infty)$ which clearly forces $\vert\nabla u\vert\notin L^p(\mathbb{R}^N\setminus\overline{B_1})$.

Suppose now that $u$ satisfies $(H_S)$. We will show that $u$ satisfies $(S')$. Since $u$ satisfies $(S)$ (again from Proposition \ref{smallproposition}), what is left to show in this case is that $\vert\nabla u\vert\in L^p(\mathbb{R}^N\setminus\overline{B_1})$ for every $N/\left(N/2+\sqrt{N-1}-1\right)<p\leq\infty$. For this purpose, let us observe that from standard regularity theory, since $u$ is bounded we have that $\vert\nabla u\vert$ is also bounded. Then, by interpolation, it suffices to prove our claim for $N/\left(N/2+\sqrt{N-1}-1\right)<p<2$. To this end consider $r\geq 1$. Applying Lemma \ref{smalllemma1} and H{\"o}lder inequality to functions $\vert u_r\vert^p$ and constant $1$ and conjugate exponents $2/p$ and $2/(2-p)$ we deduce

$$\int_r^{2r}s^{N-1}\vert u_r(s)\vert^p ds\leq (2r)^{N-1}\int_r^{2r}\vert u_r(s)\vert^p ds$$
$$\leq (2r)^{N-1}\left(\int_r^{2r}\vert u_r(s)^2 ds\right)^{p/2}\left(\int_r^{2r} ds\right)^{(2-p)/2}$$
$$\leq (2r)^{N-1}K^{\frac{p}{2}}r^{\frac{p}{2}\left( -N-2\sqrt{N-1}+3\right)}r^{\frac{2-p}{2}}=2^{N-1}K^{\frac{p}{2}}r^{\, p\left(-N/2-\sqrt{N-1}+1\right)+N},$$

\noindent for a constant $K>0$ not depending on $r\geq 1$. Applying this inequality to $r=2^j$, where $j\geq 0$ is an integer, we obtain

$$\int_1^\infty s^{N-1}\vert u_r(s)\vert^p ds=\sum_{j=0}^\infty \int_{2^j}^{2^{j+1}} s^{N-1}\vert u_r(s)\vert^p ds$$
$$\leq 2^{N-1}K^{\frac{p}{2}}\sum_{j=0}^\infty 2^{j\left(p\left(-N/2-\sqrt{N-1}+1\right)+N\right)}.$$

Finally, since $N\geq 3$ and  $N/(N/2+\sqrt{N-1}-1)<p<2$, then $p(-N/2-\sqrt{N-1}+1)+N<0$, which implies that the above series is convergent. It follows that $s^{N-1}\vert u_r(s)\vert^p\in L^1(1,\infty)$ which clearly shows that $\vert\nabla u\vert\in L^p(\mathbb{R}^N\setminus\overline{B_1})$ and the claim is proved.

Again, to finish the proof, let us observe that conditions $(L')$ and $(S')$ are clearly incompatible. \qed

\

\noindent\textbf{Proof of Corollary \ref{morse}.}
Let $N\geq 2$, $g\in C^1(\mathbb{R})$ and $u$ be a radial solution of (\ref{space}) which is stable outside a compact set (not necessarily bounded). Then, there exists $R_0>0$ such that $u$ is stable in $\mathbb{R}^N\setminus B_{R_0}$. It follows easily that the function $w(x):= u(R_0 x)$ is an stable solution of $-\Delta w=R_0^2 g(w)$ in $\mathbb{R}^N\setminus B_1$ and we can apply Theorem \ref{dichotomy} to $w$. The proof is complete by observing that there exists $u_\infty:=\lim_{r\to \infty}u(r)=\lim_{r\to \infty}w(r)$ and that $u$ satisfies $(L)$ if and only if $w$ satisfies $(L)$, while  $u$ satisfies $(S)$ if and only if $w$ satisfies $(S)$. \qed

\

\noindent\textbf{Proof of Corollary \ref{morsebounded}.}
Applying Theorem \ref{dichotomybounded}, this follows by the same method as in Corollary \ref{morse}. \qed

\section{Optimality of the main results and final remarks}

We will see that the main results obtained in this paper are optimal. To this end, for every $N\geq 2$, let us define a family $\{ u_\alpha , \alpha\in \mathbb{R} \} \subset C^\infty (\mathbb{R}^N\setminus B_1)$ of radial functions as
\begin{equation}\label{example}
\begin{array}{ll} u_\alpha (r)=r^\alpha \,  \ \
\ \ \  \forall r\geq 1\, , & \ \ \ \ \ \ \ \ \mbox{ if } \alpha
\neq 0.
\\ \\

u_0 (r)=\log r\,  \ \ \ \forall r\geq 1.

\end{array}
\end{equation}

\

It is easily seen that $u_\alpha$ is a solution of (\ref{mainequation}) with $f=f_\alpha\in C^1(\mathbb{R})$ defined by

$$
\begin{array}{ll}
f_\alpha (s)=\left\{ \begin{array}{ll}-\alpha (\alpha+N-2)s^{1-2/\alpha} \ \ \ \ \ \ \, &
\mbox{  if }s>0 \\ 0 & \mbox{  if }s\leq 0
\end{array} \right. & \mbox{ if } \alpha <0; \\ \\

f_\alpha(s)=\left\{ \begin{array}{ll} -\alpha (\alpha+N-2)s^{1-2/\alpha} &
\mbox{  if }s\geq 1 \\ (\alpha+N-2)\left( (2-\alpha)s-2\right) & \mbox{ if
}s<1
\end{array} \right. & \mbox{  if } \alpha >0; \\ \\

f_0(s)=-(N-2)e^{-2s}\, \ \ \ \ \ \ \  \ \ \ \ \ \ \ \ \ \ \ \ \  \ \ \ \mbox{ if } s\in \mathbb{R}.
\end{array}
$$

The following example shows that the exponents $-N/2\pm \sqrt{N-1}+2$ which appear in Theorem \ref{dichotomy} are optimal.

\begin{example}\label{optimalexample}

For $N\geq 2$ consider the family $\{ u_\alpha, \alpha \in \mathbb{R}\}$ defined in (\ref{example}). Then,

$$u_\alpha \mbox{ is stable }\Leftrightarrow \left( \alpha\geq -N/2+\sqrt{N-1}+2 \ \mbox{ or }\ \alpha\leq -N/2-\sqrt{N-1}+2\right) .$$

\end{example}

\noindent\textbf{Proof.}
Consider the above-mentioned functions $f_\alpha$, $\alpha\in \mathbb{R}$. We check at once that

$$f'_\alpha (u_\alpha(r))=\frac{-(\alpha-2)(\alpha+N-2)}{r^2} \ \ \ \forall r\geq 1, \, \alpha\in\mathbb{R}.$$

Consider now Hardy Inequality: $\int_{\mathbb{R}^N} ((N-2)^2/(4r^2)) v^2\leq \int_{\mathbb{R}^N} \vert \nabla v \vert^2$, for every $v \in C^1 (\mathbb{R}^N)$ with compact support, for $N\geq 3$. It is well known that the coefficient $(N-2)^2/4$ is optimal. Moreover it follows easily that this coefficient is also optimal if we consider $v \in C^1 (\mathbb{R}^N\setminus \overline{B_1})$ with compact support in $\mathbb{R}^N\setminus \overline{B_1}$, for $N\geq 2$. Hence, the stability of $u_\alpha$ is equivalent to

$$-(\alpha-2)(\alpha+N-2)\leq \frac{(N-2)^2}{4},$$

\noindent which is equivalent to

$$\alpha\geq -N/2+\sqrt{N-1}+2 \ \mbox{ or }\ \alpha\leq -N/2-\sqrt{N-1}+2.$$ \qed

\

\textbf{Remark 1.} As we have mentioned, the above example shows the optimality of Theorem \ref{dichotomy}. Moreover, taking this type of solutions, it is also possible to demonstrate the optimality of Theorem \ref{dichotomybounded}. To this purpose, consider $\alpha=-N/2+\sqrt{N-1}+2$ and $N\geq 11$. Then it is a simple matter to see that $u_\alpha$ is a bounded radial stable solution satisfying $\vert\nabla u\vert\in L^p(\mathbb{R}^N\setminus\overline{B_1})$ for every $\frac{N}{N/2-\sqrt{N-1}-1}<p\leq \infty$. This proves the optimality of Theorem \ref{dichotomybounded} for solutions $u$ satisfying $(L')$. To see the optimality of Theorem \ref{dichotomybounded} for solutions $u$ satisfying $(S')$ consider $\alpha=-N/2-\sqrt{N-1}+2$ and $N\geq 3$. Then it is easily seen that $u_\alpha$ is a bounded radial stable solution satisfying $\vert\nabla u\vert\notin L^p(\mathbb{R}^N\setminus\overline{B_1})$ for every $1\leq p\leq\frac{N}{N/2+\sqrt{N-1}-1}$.

\

\textbf{Remark 2.} Corollaries \ref{morse} and \ref{morsebounded} are also optimal. To see this, let us observe that it is possible to extend the family $\{ u_\alpha , \alpha\in \mathbb{R} \} \subset C^\infty (\mathbb{R}^N\setminus B_1)$ of radial functions defined by (\ref{example}) to another family $\{ \overline{u_\alpha} , \alpha\in \mathbb{R} \} \subset C^\infty (\mathbb{R}^N)$ of radial functions satisfying $\overline{u_\alpha}(r)=u_\alpha(r)$, for every $r\geq 1$, $\alpha\in \mathbb{R}$, such that $\overline{u_\alpha}$ is a solution of (\ref{space}) for some $g=g_\alpha\in C^1(\mathbb{R})$. To see this consider the following functions:

If $\alpha<0$, take a $C^\infty$ radial function $\overline{u_\alpha}$ satisfying $\overline{u_\alpha}(r)=2-r^2$, if $r\in [0,1/2]$, $\overline{u_\alpha}(r)=u_\alpha(r)$, if $r\in [1,\infty)$ and $\overline{u_\alpha}'(r)<0$, if $r>0$. Then  $\overline{u_\alpha}$ is a solution of (\ref{space}) for $g=g_\alpha\in C^1(\mathbb{R})$ defined by $g_\alpha(s)=2N$, if $s>2$; $g_\alpha(s)=-\Delta \overline{u} (\overline{u}^{-1}(s))$, if $0<s\leq 2$; $g_\alpha(s)=0$, if $s\leq 0$.

If $\alpha\geq 0$, take a $C^\infty$ radial function $\overline{u_\alpha}$ satisfying $\overline{u_\alpha}(r)=r^2-1$, if $r\in [0,1/2]$, $\overline{u_\alpha}(r)=u_\alpha(r)$, if $r\in [1,\infty)$ and $\overline{u_\alpha}'(r)>0$, if $r>0$. Then  $\overline{u_\alpha}$ is a solution of (\ref{space}) for $g=g_\alpha\in C^1(\mathbb{R})$ defined by $g_\alpha(s)=-2N$, if $s<-1$; $g_\alpha(s)=-\Delta \overline{u} (\overline{u}^{-1}(s))$, if $s\geq -1$.

We claim that $\overline{u_\alpha}$ is stable outside a compact set $K\subset\mathbb{R}^N$ if and only if $\alpha\geq -N/2+\sqrt{N-1}+2$ or $\alpha\leq -N/2-\sqrt{N-1}+2$. Indeed, the sufficient condition follows from Example \ref{optimalexample} (we can take $K=\overline{B_1}$) and the necessary condition is deduced from Corollary \ref{morse}.

Let us emphasize that, by a result of the author (\cite{yo}), the solutions $\overline{u_\alpha}$ are unstable in $\mathbb{R}^N$ for $\alpha\leq -N/2-\sqrt{N-1}+2$.

\

\textbf{Remark 3.} Lemma \ref{monotonicity} says that $u_r$ vanishes at most in one value in $[1,+\infty)$, for every nonconstant radial stable solution of (\ref{mainequation}), where $N\geq 2$ and $f\in C^1(\mathbb{R})$. In fact it is possible to prove that if $u$ satisfies $(S)$, then $u_r$ does not vanish  in $[1,+\infty)$. Indeed, suppose by contradiction that $u_r(r_1)=0$ for some $1\leq r_1<\infty$. Since $u$ satisfies $(H_S)$ (from Propositions \ref{largeproposition} and \ref{smallproposition}), we can take $R_1=r_1$ and apply Lemma \ref{maintool} to $\eta_0\in C^{0,1}([R_1,R_2])$, obtaining a contradiction.

\

\textbf{Remark 4.}
Let $N\geq 2$ and $p>1$. In \cite[Th. 5]{farina2} it is stated that if $0\not\equiv u\in C^2(\mathbb{R}^N)$ is a radial solution of $-\Delta u=\vert u\vert^{p-1} u$ in $\mathbb{R}^N$ which is stable outside a compact set of $\mathbb{R}^N$, only two cases occur:
\begin{enumerate}
\item[(a)] $N\geq 3$, $p=\frac{N+2}{N-2}$, $u(r)=\epsilon \left(\frac{\lambda\sqrt{N(N-2)}}{\lambda^2+r^2}\right)^{(N-2)/2}$

\noindent with $\lambda>0$, $\epsilon\in\{-1,1\}$.
\item[(b)] $N\geq 11$, $p\geq p_c:=\frac{(N-2)^2-4N+8\sqrt{N-1}}{(N-2)(N-10)}$, $u$ is stable$, u(r)=\epsilon \alpha^{2/(p-1)}v(\alpha r)$

\noindent with $\alpha>0$, $\epsilon\in\{-1,1\}$. The profile $v$ satisfies: $v(0)=1$, $v>0$, $v'<0$ in $\mathbb{R}^+$.

According to our classification of radial solutions which are stable outside a compact set (Corollary \ref{morse}) we see at once that the solutions of (a) satisfy (S), while solutions of (b) satisfy (L).

\end{enumerate}

\end{document}